\documentclass[11pt,a4paper,dvipdfmx]{amsart}
\usepackage[dvipdfm,truedimen,left=25truemm,right=25truemm,top=30truemm,bottom=30truemm]{geometry}
\usepackage{amsmath, amsfonts, amssymb, amsthm, amscd, latexsym, mathtools}
\usepackage{ascmac}
\usepackage{enumerate}
\usepackage{comment}
\usepackage{pb-diagram}
\usepackage{amsthm}
\usepackage{bm}
\usepackage{xcolor}
\usepackage{cleveref}
\usepackage{mathrsfs}
\usepackage{physics}

\numberwithin{equation}{section}

\theoremstyle{definition}
\newtheorem{Def}{Definition}[section]
\crefname{Prop}{Definition}{Definition}
\newtheorem{Prop}[Def]{Proposition}
\crefname{Prop}{Proposition}{Proposition}
\newtheorem{Thm}[Def]{Theorem}
\crefname{Thm}{Theorem}{Theorem}

\crefname{Lem}{Lemma}{Lemma}
\newtheorem{Rem}[Def]{Remark}
\crefname{Rem}{Remark}{Remark}

\crefname{Cor}{Corollary}{Corollary}

\crefname{Ex}{Example}{Example}

\crefname{Fact}{Fact}{Fact}
\newtheorem{Question}[Def]{Question}
\crefname{Question}{Question}{Question}

\newtheorem{problem}{Problem}
\crefname{problem}{Problem}{Problem}

\newcommand{\RR}{\mathbb R}

\newcommand{\Sym}{\mathrm{Sym}}

\newcommand{\T}{\mathop{\mathsf{T}}}
\newcommand{\diag}{\mathop{\mathrm{diag}}}

\newcommand{\Ad}{\mathop{\mathrm{Ad}}}
\newcommand{\Lie}{\mathop{\mathrm{Lie}}}
\newcommand{\Exp}{\mathop{\mathrm{Exp}}}

\newenvironment{customproof}[1][Proof]{%
  \par\noindent\textit{#1.} \ignorespaces
}{%
  \qed
}

\title[Invariant Conjugate Symmetric Statistical Structures on $\mathcal{N}_0$]{On Invariant Conjugate Symmetric Statistical Structures on the Space of Zero-Mean Multivariate Normal Distributions}
\author{Hikozo Kobayashi and Takayuki Okuda}
\address{Graduate School of Advanced Science and Engineering, Hiroshima University, 1-3-1 Kagamiyama, Higashi-Hiroshima City, Hiroshima, 739-8526, Japan}
\email{hikozo-kobayashi@hiroshima-u.ac.jp, okudatak@hiroshima-u.ac.jp}


\begin{document}

\begin{abstract}
By the results of Furuhata--Inoguchi--Kobayashi [Inf. Geom. (2021)] and Kobayashi--Ohno [Osaka Math. J. (2025)], the Amari--Chentsov $\alpha$-connections on the space $\mathcal{N}$ of all $n$-variate normal distributions are uniquely characterized by the invariance under the transitive action of the affine transformation group among all conjugate symmetric statistical connections with respect to the Fisher metric. 
In this paper, we investigate the Amari--Chentsov $\alpha$-connections on the submanifold $\mathcal{N}_0$ consisting of zero-mean $n$-variate normal distributions. It is known that $\mathcal{N}_0$ admits a natural transitive action of the general linear group $GL(n,\RR)$.
We establish a one-to-one correspondence between the set of $GL(n,\RR)$-invariant conjugate symmetric statistical connections on $\mathcal{N}_0$ with respect to the Fisher metric and the space of homogeneous cubic real symmetric polynomials in $n$ variables. 
As a consequence, if $n \geq 2$, we show that the Amari--Chentsov $\alpha$-connections on $\mathcal{N}_0$ are not uniquely characterized by the invariance under the $GL(n,\RR)$-action among all conjugate symmetric statistical connections with respect to the Fisher metric.
Furthermore, we show that any invariant statistical structure on a Riemannian symmetric space is necessarily conjugate symmetric.
\end{abstract}

\maketitle

\section{Introduction}\label{section:introduction}

Throughout this paper, we adopt the formulation of \emph{statistical structures} on manifolds as a pair consisting of a Riemannian metric and a symmetric $(0,3)$-tensor field (cf.~\cite{L-SM}). This formulation is equivalent to the definition as a pair of a Riemannian metric and a torsion-free affine connection compatible with it. 
A statistical structure $(g,C)$ on a smooth manifold $M$ is said to be \emph{conjugate symmetric} if the $(0,4)$-tensor field $\nabla^g C$ is symmetric (see~\cite{L-SM,MTA_2006} for details), where $\nabla^g$ denotes the Levi-Civita connection associated with $g$.
For each Riemannian manifold $(M,g)$, we denote by $S^3(T^\ast M)_{g\mathchar`-\mathrm{CS}}$ the subspace of the space $S^3(T^\ast M)$ of symmetric $(0,3)$-tensor fields on $M$ consisting of those $C$ for which the pair $(g,C)$ defines a conjugate symmetric statistical structure.

Let $M$ be an exponential family.
We denote by $g^F$ the Fisher metric on $M$ and $C^{A(\alpha)}$ the Amari--Chentsov $\alpha$-tensor field on $M$ ($\alpha \in \RR$).
Then it is well-known that the statistical structure $(g^F,C^{A(\alpha)})$ is conjugate symmetric (cf.~\cite{Amari_1985,L-SM}).

In this paper, we are concerned with the following problem:

\begin{problem}\label{problem:1}
In the setting above, 
find a characterization of the Amari--Chentsov $\alpha$-tensor fields $C^{A(\alpha)}$ on $(M,g^F)$ 
among $S^3(T^\ast M)_{g^F\mathchar`-\mathrm{CS}}$.
\end{problem}

One well-known answer to Problem \ref{problem:1} is the \emph{generalization of Chentsov's theorem} (cf.~\cite[Corollary 5.3 in Chapter 5]{Ay_2017}).
On the other hand, we focus in particular on the ``symmetry'' of the fixed space $M$, in this paper.

By Furuhata--Inoguchi--Kobayashi \cite{FIK} (for $n = 1$) and Kobayashi--Ohno \cite{KO_2025} (for $n \geq 2$), 
the Amari--Chentsov $\alpha$-tensor fields
on the $n$-variate normal distribution family
\begin{equation*}
    \mathcal{N} := \{N(x \mid \mu, \Sigma) \mid \mu \in \RR^n,~\Sigma \in \Sym^+(n,\RR)\} \cong \RR^n \times  \Sym^+(n, \RR)
\end{equation*}
are known to be characterized by $\mathrm{Aff}(n,\RR)$-invariance among $S^3(T^\ast \mathcal{N})_{g^F\mathchar`-\mathrm{CS}}$, 
where $\mathrm{Sym}^+(n, \RR)$ is the space of all positive definite symmetric matrices of order $n$, 
 and $N(x \mid \mu, \Sigma)$ denotes the $n$-variate normal distribution with mean vector $\mu$ and variance-covariance matrix $\Sigma$.
Furthermore, our previous work \cite{K-O-O-T_2025} showed that the Amari--Chentsov $\alpha$-tensor fields on the exponential family 
\begin{equation*}
    \mathcal{N}_T := \{N(x \mid \mu, \diag(\sigma^2, \dots, \sigma^2)) \in \mathcal{N} \mid \mu \in \RR^n, \sigma >0\} \cong \RR^n \times \RR_{>0}
\end{equation*}
are characterized by the invariance of the natural $\RR_{>0} \ltimes \RR^n$-action on $\mathcal{N}_T$ defined by 
\begin{equation*}
    (a, b) . N(x \mid \mu,  \diag(\sigma^2, \dots, \sigma^2) ) = N(x \mid a \mu + b,  \diag((a\sigma)^2, \dots, (a\sigma)^2)),
\end{equation*}
where $(a,b) \in \RR_{>0} \ltimes \RR^n$, among $S^3(T^\ast \mathcal{N}_T)_{g^F\mathchar`-\mathrm{CS}}$.

In this paper, we focus on the exponential family of zero-mean $n$-variate normal distributions, denoted by $\mathcal{N}_0$, defined as 
\begin{equation*}
    \mathcal{N}_0 := \{ N(x \mid 0, \Sigma) \mid \Sigma \in \mathrm{Sym}^+(n, \RR) \}.
\end{equation*}
Note that $\mathcal{N}_0$ can be identified with the parameter space $\mathrm{Sym}^+(n,\RR)$, and
$GL(n, \RR)$ acts naturally on $\mathcal{N}_0$ as below,
\begin{equation*}
    h . N(x \mid 0,\Sigma) = N(x \mid 0, h \Sigma h^{\mathsf{T}}),
\end{equation*}
where $h \in GL(n, \RR)$ and $h^{\T}$ denotes the transpose of $h$.
It is well-known that both the Fisher metric $g^F$ and the Amari–Chentsov $\alpha$-tensor field $C^{A(\alpha)}$ on $\mathcal{N}_0$ are $GL(n,\RR)$-invariant, which is a consequence of the generalization of Chentsov's theorem.

As an approach to Problem~\ref{problem:1} for $M = \mathcal{N}_0$, we examine the following question:

\begin{Question}\label{question:AC_charN_0}
Are the Amari--Chentsov $\alpha$-tensor fields on $\mathcal{N}_0$ characterized by the $GL(n,\RR)$-invariance among $S^3(T^\ast \mathcal{N}_0)_{g^F\mathchar`-\mathrm{CS}}$?
\end{Question}

The goal of this paper is to give an answer to Question \ref{question:AC_charN_0}.
The following theorem is the main theorem of this paper:

\begin{Thm}\label{thm:N0-CS}
Let $G = GL(n,\RR)$.
Let us define the vector space
\[
S^3(T^\ast \mathcal{N}_0)^{G} 
    := \{ C \in S^3(T^\ast \mathcal{N}_0) \mid C \text{ is } G\text{-invariant } \}
\]
and its linear subspace 
\[
    S^3(T^\ast \mathcal{N}_0)^{G}_{g^F\mathchar`-\mathrm{CS}} := \{ C \in S^3(T^\ast \mathcal{N}_0)_{g^F\mathchar`-\mathrm{CS}} \mid C \text{ is } G\text{-invariant } \}.
    \]
Then the following holds:
\begin{enumerate}[(1)]
    \item \label{item:main:anyCS} Any $G$-invariant statistical structure $(g,C)$ on $\mathcal{N}_0$ is conjugate symmetric. In particular, the equality $S^3(T^\ast \mathcal{N}_0)^{G}_{g^F\mathchar`-\mathrm{CS}} = S^3(T^\ast \mathcal{N}_0)^{G}$ holds. 
    \item \label{item:main:correspondence} There exists a linear isomorphism $\Phi$ from $S^3(T^\ast \mathcal{N}_0)^{G}$ onto the space $\mathcal{SP}_n^3$ of all $n$-variable  homogeneous cubic symmetric polynomials over $\RR$ such that
\[
\Phi(C^{A(\alpha)}) = \alpha (x_1^3 + \dots + x_n^3) \quad (\alpha \in \RR).
\]
    \item \label{item:main:result} The dimension of $S^3(T^\ast \mathcal{N}_0)^{G}_{g^F\mathchar`-\mathrm{CS}}$ is given as
\[
\dim S^3(T^\ast \mathcal{N}_0)^{G}_{g^F\mathchar`-\mathrm{CS}} = \begin{cases}
    3 \quad (\text{if } n \geq 3), \\
    2 \quad (\text{if } n = 2), \\
    1 \quad (\text{if } n = 1).
\end{cases}
\]
\end{enumerate}
\end{Thm}

Theorem \ref{thm:N0-CS} gives an affirmative answer to Question \ref{question:AC_charN_0} when $n=1$, and a negative one when $n \geq 2$.
We also note that a concrete example of bases of 
$S^3(T^\ast \mathcal{N}_0)^{G}_{g^F\mathchar`-\mathrm{CS}}$ can be found in Section  \ref{section:proofmain}.

We note that in the proof of Theorem \ref{thm:N0-CS} \eqref{item:main:anyCS}, we will show that 
for each symmetric space $M = G/K$, 
any $G$-invariant statistical structure $(g,C)$ on $M$ is necessarily conjugate symmetric (see Section \ref{section:symmetric} for details).
We believe that this result provides a  contribution to the study of homogeneous statistical manifolds (cf.~\cite{IO-2024}).

\section{Conjugate Symmetries on Invariant Statistical Structures on Symmetric Spaces}\label{section:symmetric}

Let $G$ be a Lie group and $M$ a homogeneous $G$-space, 
that is, $M$ is a smooth manifold equipped with a transitive smooth $G$-action.
For each point $p \in M$, 
we shall denote by $K = K^p := \{ h \in G \mid h.p = p \}$ the isotropy subgroup of $G$ at the point $p$.
Then $M$ can be regarded as the coset manifold $G/K$ via the $G$-equivariant map $G/K \rightarrow M, ~ hK \mapsto h.p$.  

The purpose of this section is to show the following theorem:

\begin{Thm}\label{theorem:cn_invmet}
In the setting above, suppose that for any (or equivalently, for some) $p \in M$, the pair $(\mathfrak{g},\mathfrak{k}^p) := (\Lie(G),\Lie(K^p))$ is a symmetric pair, 
that is, 
there exists an involutive automorphism $\theta^p$ on $\mathfrak{g}$ such that $\mathfrak{k}^p = \{ X \in \mathfrak{g} \mid \theta^p(X) = X \}$.
Then for any $G$-invariant statistical structure $(g,C)$ on $M$, $\nabla^g C \equiv 0$ holds, in particular, $(g, C)$ is conjugate symmetric.
\end{Thm}

\begin{customproof}
Theorem \ref{theorem:cn_invmet} follows directly from a combination of arguments presented in \cite[Chapters X and XI]{Kobayashi-Nomizu_II}.
For the reader's convenience, we provide a brief outline of the proof below.

Let us define, for each $X \in \mathfrak{g}$, a vector field $X^M \in \mathfrak{X}(M)$ by setting 
$$
(X^M)_q := \left. \frac{d}{dt} \right|_{t=0} \left( \exp(-tX).q \right) \in T_q M
$$
for each $q \in M$.  
It is well-known that the map $X \mapsto X^M$ defines a Lie algebra homomorphism from $\mathfrak{g}$ into the Lie algebra $\mathfrak{X}(M)$ of smooth vector fields on $M$.

For each $p \in M$, 
the canonical decomposition of $\mathfrak{g}$ with respect to $\mathfrak{k}^p := \Lie(K^{p})$ is denoted by $\mathfrak{g} = \mathfrak{k}^p + \mathfrak{p}^{p}$, i.e., we put $\mathfrak{p}^{p} := \{ X \in \mathfrak{g} \mid \theta^p(X) = -X \}$. 
Then 
$[\mathfrak{p}^p,\mathfrak{p}^p] \subset \mathfrak{k}^p$, 
$\mathfrak{p}^p$ is an $\Ad(K^p)$-stable  complement of $\mathfrak{k}^p$ in 
$\mathfrak{g}$, and the map 
\[
\mathfrak{p}^p \rightarrow T_{p} M, ~X \mapsto (X^M)_{p} = \left. \frac{d}{dt} \right|_{t=0} \left( \exp(-tX).p \right)
\]
defines a linear isomorphism.
For each tangent vector $v \in T_{p} M$, 
we write $X^v$ for the unique element in $\mathfrak{p}^p$ satisfying  $((X^v)^M)_{p} = v$.
The affine connection $\nabla^{\mathrm{cn}}$ on $M$, which is called the 
\emph{canonical connection} (cf.~\cite[Chapter X]{Kobayashi-Nomizu_II}), 
is defined by putting 
\[
\nabla^{\mathrm{cn}}_{v} C = (\mathcal{L}_{(X^v)^M}C)_{p}
\]
for each $p \in M$, each $v \in T_{p}M$ and each tensor field $C$ on $M$, 
where $\mathcal{L}_{(X^v)^M}$ denotes the Lie derivative by the vector field $(X^v)^M$.
By the definitions of $\nabla^{\mathrm{cn}}$ and $(X^v)^M$, 
one sees that $\nabla^{\mathrm{cn}}C \equiv 0$ 
for any $G$-invariant tensor field $C$ on $M$.
Furthermore, $\nabla^{\mathrm{cn}}$ is torsion-free. 
In fact, for each $p \in M$ and each $v,w \in T_pM$, 
we have $[(X^v)^M,(X^w)^M]_p = ([X^v,X^w]^M)_p = 0$ (since $[X^v,X^w] \in \mathfrak{k}^p$ and $(X^M)_p = 0$ if $X \in \mathfrak{k}^p$), 
and $\nabla^{\mathrm{cn}}_v (X^w)^M = [(X^v)^M,(X^w)^M]_p = 0$.
Hence 
\begin{align*}
(T^{\nabla^{\mathrm{cn}}})_p(v,w) 
    &= \nabla^{\mathrm{cn}}_v ((X^w)^M) - \nabla^{\mathrm{cn}}_w ((X^v)^M) - [(X^v)^M,(X^w)^M]_{p} \\
    &= 0,   
\end{align*}
where $T^{\nabla^{\mathrm{cn}}}$ denotes the torsion tensor of the affine connection $\nabla^{\mathrm{cn}}$.

Let us fix a $G$-invariant statistical structure $(g,C)$ on $M$.
Then by the invariance of the metric tensor field $g$, we have $\nabla^{\mathrm{cn}} g \equiv 0$. 
Further, $\nabla^g = \nabla^{\mathrm{cn}}$ holds since $\nabla^{\mathrm{cn}}$ is torsion-free.
By the invariance of the $(0,3)$-tensor field $C$, 
\[
\nabla^g C \equiv \nabla^{\mathrm{cn}} C \equiv 0.
\]
This completes the proof.
\end{customproof}

\section{Proof of Theorem \ref{thm:N0-CS}}\label{section:proofmain}

Let us identify $\mathcal{N}_0$ with the manifold $\Sym^+(n,\RR)$ by the correspondence $N(x \mid 0, \Sigma) \mapsto \Sigma$. 
The identity matrix of size $n$ will be denoted by $I_n \in \Sym^+(n,\RR)$.
Then $I_n$ corresponds to the standard normal distribution $N(x \mid 0, I_n)$ on $\RR^n$.
Since $\mathcal{N}_0 = \Sym^+(n,\RR)$ is an open submanifold of the vector space $\Sym(n,\RR)$, 
we have the linear isomorphism 
\[
\eta : \Sym(n,\RR) \rightarrow T_{I_n} \mathcal{N}_0,~ A \mapsto A_\eta := \left. \frac{d}{dt} \right|_{t=0} (I_n + t A).
\]

The following proposition is well-known:

\begin{Prop}[see \cite{Mitchell_1989,Skovgaard_1984}]\label{prop:FisherACalpha}
Under the identification $\eta$ above, the Fisher metric $g^F_{I_n}$ and the Amari--Chentsov $\alpha$-tensor $C^{A(\alpha)}_{I_n}$ on $T_{I_n} \mathcal{N}_0 \cong \mathrm{Sym}(n,\RR)$
can be written as below:
\begin{align}
    g^F_{I_n} &: \mathrm{Sym}(n,\RR) \times \mathrm{Sym}(n,\RR) \rightarrow \RR, ~ (X,Y) \mapsto \frac{1}{2}\tr(XY),\\
    C^{A(\alpha)}_{I_n} &: \mathrm{Sym}(n,\RR) \times \mathrm{Sym}(n,\RR) \times \mathrm{Sym}(n,\RR) \rightarrow \RR, ~ (X,Y,Z) \mapsto \alpha \cdot \tr(XYZ). \label{eq:ACalpha}
\end{align}    
\end{Prop}

Let us give a proof of Theorem~\ref{thm:N0-CS} as below:

\begin{customproof}[Proof of Theorem~\ref{thm:N0-CS}]
We put $G = GL(n,\RR)$.
Recall that $\mathcal{N}_0 = \Sym^+(n,\RR)$ is a homogeneous $G$-space equipped with the action defined by 
\[
h.\Sigma := h \Sigma h^{\T} \quad (\text{for } h \in G = GL(n,\RR),~\Sigma \in \Sym^+(n,\RR)).
\]
The isotropy subgroup $K = K^{I_n}$ of $G$ at the point $I_n \in \Sym^+(n,\RR)$ is the orthogonal group $O(n)$.
It is well-known that $(G,K) = (GL(n,\RR),O(n))$ is a symmetric pair of Lie groups, 
and hence $(\mathfrak{g},\mathfrak{k}) := (\Lie(G),\Lie(K))$ is also a symmetric pair.
Thus the claim \eqref{item:main:anyCS} in Theorem~\ref{thm:N0-CS} is followed immediately by Theorem \ref{theorem:cn_invmet}.

Let us give a proof of the claim \eqref{item:main:correspondence} in Theorem~\ref{thm:N0-CS}.
The natural $K$-action on the tangent space $T_{I_n} \mathcal{N}_0$ is called the isotropy representation at the point $I_n$.
We write $S^3(T^\ast_{I_n} \mathcal{N}_0)^K$ for the space of $K$-invariant symmetric $3$-tensors on the cotangent space at $I_n$, i.e., on $T^\ast_{I_n} \mathcal{N}_0$.
Then by the general theory of invariant sections of equivariant vector bundles over homogeneous spaces, 
one sees that the map 
\begin{equation}
S^3(T^\ast \mathcal{N}_0)^G \rightarrow S^3(T^\ast_{I_n} \mathcal{N}_0)^K,~ C \mapsto C_{I_n} \label{eq:Fiber}    
\end{equation}
gives a linear isomorphism.    

We shall define the $K = O(n)$-representation on the vector space $\Sym(n,\RR)$ by putting
\[
k . X := k X k^{-1} \quad (\text{for } X \in \Sym(n,\RR), ~ k \in O(n)).
\]
The vector space of all $K$-invariant symmetric $3$-tensors on the space $\Sym(n,\RR)^\ast$ is denoted by $S^3(\Sym(n,\RR)^\ast)^K$.
One can easily see that the identification $\eta : \Sym(n,\RR) \rightarrow T_{I_n} \mathcal{N}_0$ is an isomorphism between $K = O(n)$-representations, 
where $T_{I_n} \mathcal{N}_0$ is considered as the isotropy representation of $K$.
Therefore, $S^3(T_{I_n}^\ast M)^K$ can be identified with $S^3(\Sym(n,\RR)^\ast)^K$.
By combining this, the isomorphism \eqref{eq:Fiber} above, and Proposition  \ref{prop:FisherACalpha} \eqref{eq:ACalpha}, 
we have a linear isomorphism from $S^3(T^\ast M)^G$ onto $S^3(\Sym(n,\RR)^\ast)^K$ such that $C^{A(\alpha)}$ maps to 
the tensor $C^\alpha$ defined by 
\[
C^\alpha(X,Y,Z) := \alpha \cdot \tr(XYZ) \quad (X,Y,Z \in \Sym(n,\RR)).
\]
To complete the proof of the claim \eqref{item:main:correspondence}, 
we only need to find a linear isomorphism from $S^3(\Sym(n,\RR)^\ast)^K$ onto $\mathcal{SP}^3_n$ 
such that $C^\alpha$ maps to the polynomial $\alpha \cdot (\sum_{i} x_i^3)$.
For each $C \in S^3(\Sym(n,\RR)^\ast)^K$, 
we define the polynomial function 
$q_C$ on the vector space $\Sym(n,\RR)$ by 
\[
q_C : \Sym(n,\RR) \rightarrow \RR, ~ X \mapsto C(X,X,X).
\]
The correspondence $C \mapsto q_C$ gives a linear isomorphism between $S^3(\Sym(n,\RR)^\ast)^K$ 
and the vector space $\mathcal{P}^3(\Sym(n,\RR))^K$ of $K$-invariant homogeneous cubic polynomial functions on $\Sym(n,\RR)$.
Note that $C^\alpha$ maps to the function 
\[
q_{\alpha} : \Sym(n,\RR) \rightarrow \RR, ~ X \mapsto \alpha \cdot \tr(X^3).
\]
Furthermore, we shall write $D$ for the linear subspace of $\Sym(n,\RR)$ consisting of all diagonal matrices.
Then the symmetric group $\mathfrak{S}_n$ of order $n$ acts on $D$ by permutations of subscripts.
Let us denote by $\mathcal{P}^3(D)^{\mathfrak{S}_n}$ the space of all $\mathfrak{S}_n$-invariant homogeneous cubic polynomial functions on the vector space $D$.
We shall consider the linear isomorphism 
\[
\mathcal{SP}^3_n \rightarrow \mathcal{P}^3(D)^{\mathfrak{S}_n}, ~ P(x_1,\dots,x_n) \mapsto f_P, 
\]
where the function $f_P$ is defined by 
\[
f_P : D \rightarrow \RR, ~ \diag(\lambda_1,\dots,\lambda_n) \mapsto P(\lambda_1,\dots,\lambda_n).
\]
Then $\mathcal{SP}^3_n$ can be identified with the space $\mathcal{P}^3(D)^{\mathfrak{S}_n}$.
Note that the polynomial $\alpha \cdot (\sum_i x_i^3)$ corresponds to the function 
\[
f_{\alpha} : D \rightarrow \RR, ~ \diag(\lambda_1,\dots,\lambda_n) \mapsto \alpha \cdot \sum_{i} \lambda_i^3.
\]
Thus our goal is to find a linear isomorphism $\varphi$ from $\mathcal{P}^3(\Sym(n,\RR))^K$ onto $\mathcal{P}^3(D)^{\mathfrak{S}_n}$ such that $\varphi(q_{\alpha}) = f_{\alpha}$.
For each function $q \in \mathcal{P}^3(\Sym(n,\RR))^K$, define $\varphi(q) := q|_{D}$ by the restriction of $q$ on the linear subspace $D$. 
One sees that the correspondence $q \mapsto \varphi(q)$ defines a linear map $\varphi$ from
$\mathcal{P}^3(\Sym(n,\RR))^K$ to $\mathcal{P}^3(D)^{\mathfrak{S}_n}$, and $\varphi(q_{\alpha}) = f_{\alpha}$.
Furthermore, the map $\varphi$ is injective since for any $X \in \Sym(n,\RR)$, there exists $k \in K$ such that $\Ad(k) X \in D$ (i.e.,~any symmetric matrix is diagonalizable by an orthogonal matrix).
Therefore, it is enough to show that the map $\varphi$ is surjective.
Let us define the symmetric polynomial function $p_k$ on $D$ ($k \in \mathbb{Z}_{\geq 0}$) by 
\begin{align*}
    p_k(\diag(\lambda_1,\dots,\lambda_n)) := \sum_i \lambda_i^k \quad (\text{for } \diag(\lambda_1,\dots,\lambda_n) \in D).
\end{align*}
Then by the theory of symmetric polynomials, 
one can check that 
$\mathcal{P}^3(D)^{\mathfrak{S}_n}$
is spanned by the three homogeneous cubic polynomial functions $p_3$, $p_2 p_1$ and  $p_1^3$.
On the other hand, let us define the $K$-invariant homogeneous cubic polynomial functions $q_1,q_2,q_3$ on $\Sym(n,\RR)$ by 
\begin{equation*}
    q_1(X) := \tr(X^3), ~
    q_2(X) := \tr(X^2) \cdot \tr(X), ~
    q_3(X) := (\tr(X))^3.
\end{equation*}
Then $\varphi(q_1) = p_3$, $\varphi(q_2) = p_2 p_1$ and $\varphi(q_3) = p_1^3$.
This completes the proof of the claim \eqref{item:main:correspondence}.

The claim \eqref{item:main:result} follows from the claims \eqref{item:main:anyCS}, \eqref{item:main:correspondence} and the well-known fact that the vector space $\mathcal{SP}_n^3$ is $3$-dimensional if $n \geq 3$, $2$-dimensional if $n = 2$, and $1$-dimensional if $n=1$.
\end{customproof}

\begin{Rem}
The following three elements form a generating set of the vector space $S^3(\Sym(n,\RR)^\ast)^K \cong S^3(T^\ast \mathcal{N}_0)^G_{g^F\mathchar`-\mathrm{CS}}$:
\begin{itemize}
\item $C_1(X,Y,Z) := \tr(XYZ)$,
\item $C_2(X,Y,Z) := (1/3) (\tr(X)\tr(YZ)+\tr(Y)\tr(XZ)+\tr(Z)\tr(XY))$,
\item $C_3(X,Y,Z) := \tr(X)\tr(Y)\tr(Z)$.
\end{itemize}
In particular, if $n \geq 3$, the subset $\{C_1, C_2, C_3\}$ is a basis of the vector space $S^3(\Sym(n,\RR)^\ast)^K$.
The symmetric tensor $C_1$ corresponds to the Amari--Chentsov $+1$-tensor field on $\mathcal{N}_0$.
\end{Rem}

\begin{Rem}
It is worth emphasizing that the linear isomorphism
$\eta : \Sym(n,\RR) \rightarrow T_{I_n} \mathcal{N}_0$ 
differs from the following ``natural'' linear isomorphism:
\[
\phi : \Sym(n,\RR) \rightarrow T_{I_n} \mathcal{N}_0,\quad
A \mapsto \left. \frac{d}{dt} \right|_{t=0} \left( \Exp(-tA).I_n \right).
\]
Indeed, it can be directly verified that $\phi = -2\eta$.
\end{Rem}

\subsubsection*{Acknowledgements:} 
The authors would like to give heartfelt thanks to 
Hideyuki Ishi  
whose suggestions were of inestimable value for this paper.
The authors would also like to thank to
Hitoshi Furuhata,
Kento Ogawa,
Yu Ohno,
Hiroshi Tamaru and
Koichi Tojo
whose comments made enormous contribution to this paper.
The second author is supported by JSPS Grants-in-Aid for Scientific Research JP20K03589, JP20K14310, JP22H01124, and JP24K06714.

\end{document}